\documentstyle{amsppt}
\magnification =1200
	\hcorrection{0.25in}
\document

\topmatter 

\title  {Projective structures with discrete holonomy 
representations} 
\endtitle 

\rightheadtext{Projective structures with discrete 
holonomies} 

\author {Hiroshige Shiga and Harumi Tanigawa} 
\endauthor  

\address{Shiga: Department of Mathematics, Tokyo Institute 
of Technology, Tokyo
152 Japan}
\endaddress
\address{Tanigawa: Graduate School of Polymathematics,  Nagoya 
University,  Nagoya \hbox{464-01}
Japan} \endaddress 

\abstract  Let $K(X)$ denote the set of projective 
structures on a 
compact Riemann surface $X$ whose holonomy representations are discrete. 
We will show that each  component of the
interior of $K(X)$ is holomorphically equivalent to a complex
submanifold of the product of Teichm\"uller spaces and the holonomy
representation of every projective structure in the interior of $K(X)$ is
a quasifuchsian group.
\endabstract

\endsubjclass
\thanks {Research at MSRI is supported by NSF grant 
\#DMS--9022140}
\endthanks

\endtopmatter 
\heading
 1. Introduction
\endheading Let $\Sigma_g$ be a compact oriented 
differentiable surface
of genus
$g \ge 2$. A projective structure (or $\text {\bf 
CP}^1$-structure) on 
$\Sigma_g$ is a maximal system of charts with 
transition maps in the
projective automorphism group of $\text {\bf CP}^1$, 
namely,
$\operatorname {PSL}_2(\operatorname {\bold C})$. Let $P_g$ denote the
space of  all projective
structures on $\Sigma_g$ with markings. Since elements 
of
$\operatorname {PSL}_2(\operatorname {\bold C})$ are holomorphic, a
projective  structure determines
its underlying complex structure. Hence there is a 
natural map $\pi : P_g
\to T_g$, where $T_g$ denotes the Teichm\"uller space. 
For each
projective structure, we can take its developing map 
on the 
universal  
and have a
holomorphic quadratic differential (or a projective 
connection) on the
Riemann surface (the complex structure under the projective structure) 
 by taking the Schwarzian derivative. It is well-known (cf.
section 2 for detail) that this procedure  gives a natural
identification between
$\pi : P_g \to T_g$ and the bundle of holomorphic 
quadratic differentials
over the Teichm\"uller space $\pi : Q\to T_g$. \par

In this paper, we will investigate the set of 
projective structures on
$\Sigma_g$ with discrete holonomy representations with respect to the
parametrization by $Q$. Let 
$K \subset Q$ be
the set of all projective structures with discrete 
holonomy
representations. It is well-known that there is an 
open set $QF$ of $Q$,
which is a neighborhood of  the $0$-section $\{0\}\times 
T_g$ of $Q \to T_g$,
consisting of all projective structures with univalent 
developing maps and
quasifuchsian holonomy representations. For any 
Riemann surface $X \in
T_g$, the fiber of $\pi |QF: QF \to T_g$ over $X$, 
which we denote by
$QF(X)$,  is a domain which is holomorphically 
equivalent to the
Teichm\"uller space and is called a Bers slice. One of 
our main result
shows that a similar phenomenon occurs  for each 
component of the interior point set 
$\operatorname{int} K$. 
\par 

Let $Q(X)$ denote the fiber over $X$ of $Q \to T_g$, namely, $Q(X)$
is the space of projective structures on the complex structure $X$.
 Let
$K(X)$ denote the  fiber of $\pi |K : K
\to X$. We will denote by $\operatorname {int}_{Q(X)} K(X)$
the interior of $K(X)$ considered in $Q(X)$. We employ this notation
to avoid mixing up $\operatorname {int}_{Q(X)} K(X)$ with the fiber of
$\operatorname {int} K $ (interior is considered in $Q$) over $X$,
although it will turn out these two are the same by our result.
As we recalled above, 
$\operatorname {int}_{Q(X)} K(X)$ has
a component containing the Bers slice $QF(X)$. 
(In fact, the Bers slice
$QF(X)$  coincides with the
component of
$\operatorname {int}_{Q(X)} K(X)$ containing 0 (\cite {Sh}).) On the
other hand, it was shown by Maskit [M] that there are complex structures
$X$ on  which there exist
components of
$\operatorname {int}_{Q(X)} K(X)$  other than the Bers slice (see also
Hejhal \cite {H} and Goldman \cite {G}).  Indeed, he showed
that  given a Fuchsian group isomorphic to $\pi_1 
\Sigma_g$ there are projective structures outside the Bers slice
on {\it some} complex structures whose holonomy representation
are conjugate to the given Fuchsian group in
$\operatorname{PSL}_2(\bold {C})$. 
However, the Fuchsian holonomy representation hardly gives 
information about
the  complex structures under the projective structures 
when the projective structure is not in the Bers slice, because the
developing map is  not a covering map (\cite {Kr1}).  Here arise natural
questions: 
\roster
\item
{For {\it every} complex $X
\in T_g$ are there any components of
$\operatorname {int}_{Q(X)} K(X)$ other than the Bers slice?}
\item
      {What kind of discrete groups appear as 
holonomy representations
in such components? }
\item {What kind of analytic properties 
do such components
have? }
\endroster 

As for the first question, an affirmative answer 
is
given in  [T].  The aim of this paper is to work on 
the other two
questions. We will show that each component of
$\operatorname {int}_{Q(X)} K(X)$ consists of quasifuchsian 
groups and is
holomorphically equivalent to a submanifold of the 
product $T_g
\times T_g$.\par
The authors would like to thank Curt McMullen for stimulating
conversations and comments. Considerable part of this work was done at
Mathematical Sciences Research Institute. The authors are very grateful
for their hospitality.
\remark {Remark} All results in this paper holds for {\it bounded
projective structures  on Riemann surfaces of finite type} (i.e.
projective structures determined by bounded holomorphic quadratic
differentials) by parallel arguments.
\endremark
\head 
 2. Notation and basic facts
\endhead
In this section, we recall some known facts we will use to prove
 our results.

 Throughout this paper, $\Gamma$ denotes a 
Fuchsian group
acting on the upper half plane {\bf H} and the lower 
half plane $L$  such
that the quotient   $X = L/ \Gamma$ is a compact 
Riemann surface of genus
$g \ge 2$.

\subhead 2.1 Quasiconformal deformations and the Teichm\"uller space
\endsubhead

For a finitely 
generated Kleinian
group $G$, we denote by $\Omega (G)$ its region of 
discontinuity and by
$\Lambda (G)$ its limit set. 

\par

For a $G$-invariant open 
set $U \subset
\Omega (G)$,  a measurable $(-1,1)$-form  $\mu$ on the Riemann sphere 
with the following
properties  is called a {\it Beltrami differentials for } $G$
{\it supported on } $U$:
$$\mu|(\hat {\text {\bf C}}-U )= 0, \qquad ||\mu||_{\infty} 
< 1$$
 and  
$$\mu \circ g \times \overline {g'} /g'  =\mu  \qquad \text {a.e.}$$
 for all $g
\in G$. We denote the  set of all Beltrami 
differentials on $U$ by
$B_1(G, U)$. When
$U = \Omega (G)$, we abbreviate $B_1(G, U)$ by 
$B_1(G)$. The space
$B_1(G)$ is the unit ball of the complex Banach space 
of all  measurable
$(-1,1)$-form for $G$ with 
supremum norm. \par
By Ahlfors-Bers' theory, for each
Beltrami differential $\mu$ on $\hat {\operatorname {\bold C}}$ with
$||\mu||_\infty$, there exists a  quasiconformal mapping
$f^{\mu}$ with Beltrami differential $\mu$, that is,
$$(f^{\mu})_{\overline z} = \mu \times (f^{\mu})_z.$$
 Such a 
mapping is unique up to
post-composition of elements of
$\operatorname{PSL}_2 (\bold {C})$. In each 
argument below, we will
employ  some convenient normalization. \par
It is easy to see that if $\mu$ is in $B_1(G)$ then  
$f^{\mu} \circ g \circ (f^{\mu})^{-1}$ is a M\"obius transformation for
every $g \in G$ and $f^{\mu} G (f^{\mu})^{-1}$ is a subgroup
of  $\operatorname{PSL}_2(\bold {C})$ with region of discontinuity
$f^\mu (\Omega (G))$. \par

Two Beltrami  differentials 
$\mu_1$ and $\mu_2$ for $G$ are said to be equivalent if 
$f^{\mu_1}$ and
$f^{\mu_2}$ determine an equivalent homomorphism of 
$G$, namely, there exists a M\"obius transformation $A \in
\operatorname{PSL}_2(\bold {C})$ such that 

$$f^{\mu_1} \circ g \circ (f^{\mu_1})^{-1} = A \circ f^{\mu_2} \circ g  
\circ (f^{\mu_2})^{-1} \circ A^{-1}$$
for every $g \in G$. 
 A Beltrami differential is called {\it trivial} if it 
is equivalent to
 $0$. 
\definition {Definition 2.1}
The set of all equivalence 
classes of Beltrami
differentials for $G$ is  called the {\it quasiconformal 
deformation space of }
$G$ and is denoted by $QC(G)$. \par 

For a Fuchsian group $\Gamma$ acting on the lower half plane $L$
the {\it Teichm\"uller space} $T(\Gamma)$ of 
$\Gamma$ is the quotient space of 
$B_1(\Gamma, L)$ by
the equivalence relation as above.
\enddefinition 
The space $QC(G)$ has a 
natural complex
structure such that the canonical projection $B_1(G) 
\to QC(G)$ is
holomorphic. Also,
 the Teichm\"uller space $T(\Gamma)$ has a  natural 
complex structure
such that the quotient map $B_1(\Gamma, L) \to 
T(\Gamma)$ is holomorphic. The
Teichm\"uller space $T(\Gamma)$ is also regarded as a 
deformation space
of the Riemann surface $X =L/\Gamma$, which we
denote  by $T(X)$. It is
well-known that the complex structure of the 
Teichm\"uller space of
Riemann surfaces of genus $g$ is independent of the 
choice of the base
point $X$.
 When we need not mention the base point
$X$ we denote the Teichm\"uller space of Riemann 
surfaces of genus $g$ by
$T_g$.  \par

Here, we recall another definition 
of $T(\Gamma)$ as the space of Fuchsian groups.
For $\mu \in B_1(\Gamma, L)$, let  $f_\mu$ denote the 
quasiconformal
homeomorphism of $L$ onto itself with Beltrami differential
$\mu|L$ fixing $0,1$ and
$\infty$. The mapping $f_\mu$ is the restriction of a quasiconformal
mapping $f^{\mu'}$ where $\mu'$ is a Beltrami differential on
$\hat {\operatorname {\bold C}}$ defined by extending $\mu |L$ to 
$\hat {\operatorname {\bold C}}$
symmetrically beyond $\operatorname {\bold R}$:

$$\mu'(z)  = \cases \mu (z),  & \text {if  $z \in L$,}\\
              \overline {\mu (\overline z)}, & \text{if $ z \in H$}.
\endcases $$

Then the group $f_\mu \Gamma (f_\mu)^{-1}$ is a Fuchsian group acting on
$L$. It is easy to see that two Beltrami differentials 
$\mu_1$ and $\mu_2$
in
$B_1 (\Gamma, L)$ determines the same homomorphism $\Gamma \to
\operatorname {PSL}_2(\bold R)$ if and only if $\mu_1'$ and $\mu_2'$ are
equivalent in $B_1(\Gamma)$, where $\mu_i'$ ($i = 1,2$) is the symmetric
extension of $\mu_i$ as above. It is also easy to see that $\mu_1'$ and
$\mu_2'$ are equivalent in $B_1(\Gamma)$ if and only if $\mu_1$ and
$\mu_2$ are equivalent in $B_1(\Gamma)$.

 \definition {Definition 2.2 (Teichm\"uller space as the space of
Fuchsian groups)} The Teichm\"uller space
$T(\Gamma)$ is the equivalence classes of the symmetric Beltrami
differentials, namely, equivalence classes of Beltrami differentials $\mu'
\in B_1(\Gamma)$ such that $\mu'(z) = \overline {\mu' (\overline z)}$
almost everywhere.
\enddefinition

By the above remark,   this definition of $T(\Gamma)$  is the same as 
that in Definition 2.1, if we regard $T(\Gamma)$ as a topological space.
Actually, these two definition give the same real analytic structures.
However, when we regard $T(\Gamma)$ as the space of Fuchsian groups, it is
not considered as a complex manifold, since the
space of symmetric Beltrami differentials is not a complex Banach
manifold.\par
 
 See Lehto [L], for example, for more detail on fundamental facts on
Teichm\"uller spaces.

\par

\subhead 
2.2 Projective structures and quadratic differentials
\endsubhead
 Assume that we are given a projective structure on
$\Sigma_g$. Since the transition mappings are holomorphic, as they are
restrictions of elements of 
$\operatorname{PSL}_2(\bold {C})$, we have a complex structure under the
projective structure. Let $X$ denote the Riemann surface and let $\Gamma$
be a Fuchsian group acting on $L$ such that $X = L/\Gamma$.
Take a coordinate function of the projective structure. We can take its
analytic continuation along any curve on $X$ and have a multivalued
holomorphic mapping into the Riemann sphere. This multivalued mapping is
lifted to a locally univalent meromorphic function $W : L \to
\operatorname {\bold CP}^1$ on the universal covering space $L$. This
mapping  is called a {\it developing map} of the projective structure.
Note that the  developing map  is determined by the projective structure
uniquely up to post  compositions  of
elements of $\operatorname{PSL}_2(\bold {C})$. \par
 When we take the 
analytic continuation of  a local
coordinate function along a closed curve and  come 
back to the starting
point, the values differ  from each 
other by elements of
$\operatorname{PSL}_2(\bold {C})$, since the 
transition mappings are in
$\operatorname{PSL}_2(\bold {C})$. Therefore, we have  
a homomorphism of 
$\pi_1 X$ to $\operatorname {PSL}_2(\operatorname {\bold C})$. 
 If we look at this on the universal covering space $L$, we have a
homomorphism
$\chi :\Gamma
\to 
\operatorname {PSL}_2(\operatorname {\bold C})$ such that

$$ W \circ \gamma (z) = \chi (\gamma)\circ W(z), \quad z \in L \tag 2.1 $$
for all $\gamma \in \Gamma$.
 This 
homomorphism is called the {\it
holonomy representation}.\par
Thus the projective structure determines the pair $(W, \chi)$ uniquely up
to the action of $\operatorname {PSL}_2(\operatorname {\bold C})$.
Conversely, it is easy to see that given a pair $(W, \chi)$ of locally
univalent meromorphic function  $W : L \to \operatorname {\bold CP}^1$ and
a homomorphism $\chi : \Gamma \to \operatorname {PSL}_2(\operatorname
{\bold C})$, where
$\Gamma$ is a Fuchsian group, satisfying $(2.1)$ we have a projective
structure on the Riemann surface $L / \Gamma$. Therefore, the space of
all projective structures on $\Sigma_g$ is identified with the set of
pairs $(W, \chi)$ satisfying $(2.1)$ modulo the action of 
$\operatorname {PSL}_2(\operatorname {\bold C})$.\par

Now we recall the parametrization of the space of projective structures
by holomorphic quadratic differentials.  Given a projective structure on
a Riemann surface $X = L/\Gamma$, take the Schwarzian derivative of the
developing map $W$ and denote it by $\varphi$. (Here, Schwarzian derivative
of a locally univalent meromorphic function
$f$ is  defined by
$(f''/f')'-1/2(f''/f')^2$.) Then   
$\varphi$ is a holomorphic function on $L$ satisfying
$$\varphi = \varphi\circ\gamma \times {\gamma '}^2 \tag 2.2$$ for 
all $\gamma\in
\Gamma$. Such a holomorphic function  
satisfying $(2.2)$ 
  is called a {\it holomorphic quadratic 
differential for} $\Gamma$. A holomorphic quadratic differential for
$\Gamma$ is the lift of a holomorphic quadratic differential on the
Riemann surface $X = L/\Gamma$. We have seen that given a projective
structure on
$X$  we have a
holomorphic quadratic differential for $\Gamma$. 
Conversely, given a
holomorphic quadratic differential $\varphi$ for $\Gamma$, it is well-known
(see {\cite H} or  \cite {Kr2}) that there is a locally univalent
meromorphic function
$W_\varphi$ whose Schwarzian derivative is $\varphi$. It is easy to see
from
$(2.1)$ that there is a homomorphism $\chi_\varphi:\Gamma \to \operatorname
{PSL}_2 (\bold C)$ satisfying $(2.1)$.
 Thus the space of all projective
structure on the surface
$\Sigma_g$  is  identified with
the fiber space $Q \to T_g$ of holomorphic quadratic 
differentials on compact
Riemann surfaces of genus $g$. 
\par
\definition {Notation 2.3}  In this paper, we will denote by 
$W_\varphi$ the developing map and by $\chi_\varphi$ 
the holonomy
representation of the projective structure determined 
by $\varphi$. Since
they are determined up to elements of 
$\operatorname{PSL}_2(\bold {C})$,
we have to fix normalization in each argument in the 
following section.
The choice is not essential.
 We shall employ convenient ones each time.
\enddefinition

\definition {Notation 2.4}
 Let $Rep$ denote the space of 
$\operatorname{PSL}_2(\bold {C})$
representations of the fundamental group $\pi_1 X $ : 
$Rep  =  Hom(\pi_1X,\operatorname{PSL}_2(\bold {C})) 
/ \operatorname{PSL}_2(\bold {C})$.
 Here, two $\operatorname{PSL}_2(\bold {C})$ 
representations are
equivalent if they differ by conjugation  in 
$\operatorname{PSL}_2(\bold
{C})$. We define a  mapping $h : Q \to Rep$ by sending 
each $\varphi
\in Q$ to the holonomy representation of the 
projective structure
determined by $\varphi$, and call it the {\it holonomy map}.
\enddefinition
We will use the following property of $h$:
\proclaim {Theorem 2.5 (Hejhal \cite {H, Theorem 1})}
The mapping $h : Q \to Rep$ is not a covering map. 
However, it is  a
local $C^1$-diffeomorphism.
\endproclaim
 Note that 
when we regard $Q$ as the space of quadratic 
differentials for
Fuchsian groups on $L$, we regard $T(\Gamma)$ as the 
space of
Fuchsian groups,  hence not a complex manifold.  
In particular, the vector
bundle $Q\to T(\Gamma)$ is not holomorphic. However, 
each fiber
$Q(\Gamma')$ over
$\Gamma' \in T(\Gamma)$ has a natural complex 
structure as a finite
dimensional complex linear space.

\proclaim {Lemma 2.6} The holonomy map 
$h:Q(\Gamma) \to Rep $  is holomorphic for each $\Gamma$.
\endproclaim 

Note that the fiber $Q(\Gamma)$ is
identified with the space of projective structure whose underlying complex
structure is $X = L / \Gamma$.\par

In the next section, we will discuss 
the space of projective
structures with discrete holonomy representations.  We 
denote by $K$
the subset of $Q$ consisting of all elements $\varphi$ 
which determine
projective structures with discrete holonomy 
representations and
denote by $K(\Gamma)$ the fiber over a point $\Gamma$. We denote by
$\operatorname {int}_{Q(\Gamma)} K(\Gamma)$ the interior of $K(\Gamma)$
considered in $Q(\Gamma)$. \par
The set of all
equivalence classes of $ B_1(\Gamma, H)$ is called the {\it Bers slice}.
From the definition the Bers slice is isomorphic to the Teichm\"uller
space. For each element $\mu\in B_1(\Gamma, H)$,   $f^\mu|L$ is a
univalent meromorphic mapping and 
$f^\mu \Gamma (f^\mu)^{-1}$ is a quasifuchsian group. Therefore, it
determines a projective structure on $L/\Gamma$. Thus the Bers slice is
embedded in $Q(\Gamma)$. We will also call this image Bers slice.

 As we noted in the introduction, for each
$\Gamma$, $\operatorname {int}_{Q(\Gamma)} K(\Gamma)$ has a component which
contains the Bers slice.
\proclaim {Proposition 2.7 (\cite {Sh})} The component of 
$\operatorname {int}_{Q(\Gamma)} K(\Gamma)$ containing $0$ coincides with
the Bers slice. For any component $\kappa (\Gamma)$ 
of $\operatorname {int}_{Q(\Gamma)} K(\Gamma)$ which is not the Bers slice,
the developing map $W_\varphi$ is surjective onto the Riemann sphere for
every
$\varphi
\in \kappa (\Gamma)$. For any two elements $\varphi_1$ and $\varphi_2$ of
$\kappa (\Gamma)$, $\chi_{\varphi_1}(\Gamma)$ and $\chi_{\varphi_2}
(\Gamma)$ are quasi conformally equivalent to each other.
 
\endproclaim

\subhead 2.3 Holomorphic motions of Kleinian groups
\endsubhead 
Here we state the theorem of holomorphic motions in a convenient form
 for later use. (See [EKK] for detail 
and references for 
the theory of holomorphic motions.)
\definition {Definition 2.8} Let $G$ be a finitely generated discrete
group of $\operatorname{PSL}_2(\bold {C})$. A family  $\{\theta_z\}_{z \in
\Delta}$ of isomorphisms of $G$ into $\operatorname{PSL}_2(\bold {C})$,
parametrized by the unit disc $\Delta = \{z \in \operatorname {\bold C}
; |z|<1\}$,  is called a {\it holomorphic motion of}
$G$ or a {\it holomorphic family of isomorphisms of} $G$ if it
satisfies the following conditions:
\roster
\item $\theta_0$ is the identity: $\theta_0 (g) = g$ for all $g \in G$.
\item For every $z \in \Delta$, $\theta_z (G)$ is discrete.
\item For every $g \in G$, $z \mapsto \theta_z(g)$ is a holomorphic
mapping of $\Delta$ to  $\operatorname{PSL}_2(\bold {C})$.
\endroster 
\enddefinition
 
\proclaim  {Theorem 2.9 (Bers [B], Earle-Kra-Krushkal[EKK], Slodokowski
[Sl])} Let $\{\theta_z\}_{z \in \Delta}$ is a holomorphic family of
 isomorphisms of a 
discrete group $G$. Then there is a
holomorphic mapping $z \mapsto \mu_z$ from $\Delta$ to $B_1(G)$  
such that
$\theta_z(\gamma) = f^{\mu_z}\circ \gamma \circ 
(f^{\mu_z})^{-1}$ for all
$\gamma \in G$ where $f^{\mu_z}$ is a 
quasiconformal mapping with
Beltrami differential $\mu_z$.
\endproclaim 
\remark {Remark}
If we have a holomorphic family of  {\it discrete} groups, it is actually a
holomorphic family of {\it discontinuous} groups, unless the family is
trivial. In fact, by the above theorem the holomorphic family of
isomorphisms are given by a family of quasiconformal deformations. If the
discrete group is not discontinuous, there is no non-trivial quasiconformal
deformations by the Sullivan's rigidity theorem (\cite {Su}). Therefore,
$\operatorname {int}_{Q(\Gamma)} K(\Gamma)$ coincides with the interior of 
set of projective structures with discontinuous holonomy representations.
\endremark

\head  3. Projective structures with discrete
holonomy representations 
\endhead 

To show the main result, we will construct a local 
inverse to the
holonomy map $h : Q \to Rep$ on each component 
of $\operatorname{int}
K$, which is holomorphic in the direction of fibers of 
$Q\to T_g$.  For
the construction, we first note the following fact.

\proclaim {Lemma 3.1} Each component $\kappa (\Gamma)$ 
of
$\operatorname{int}_{Q(\Gamma)} K(\Gamma)$ is either consists of 
totally degenerate
groups or consists of quasifuchsian groups which is 
isomorphic to
$\Gamma$.
\endproclaim 
\demo {Proof} 
 
First, note that for every $\varphi \in \kappa  (\Gamma)$
$\chi_\varphi$ is an isomorphism which preserves parabolicity 
(Kra \cite {Kr2} and \cite {Kr3}). 

By the result of  Maskit \cite{M , Theorem 2 and Theorem
3}, 
$\chi_\varphi (\Gamma)$ is
either a quasifuchsian group or a totally degenerate 
group. By Proposition 2.7, either $\chi_\varphi (\Gamma)$ is
a quasifuchsian group for every $\varphi \in \kappa (\Gamma)$ or 
$\chi_\varphi (\Gamma)$ is
a totally degenerate group for every $\varphi \in \kappa (\Gamma)$.

\qed 

\enddemo

We now define a mapping which plays a key role to 
analyze the structure of
$\operatorname{int} K$ or $\operatorname{int}_{Q(\Gamma)} 
K(\Gamma)$. \par  Fix an
element $\varphi \in
\kappa (\Gamma)$ and denote the holonomy 
representation 
$\chi_\varphi (\Gamma)$ by $G_\varphi$, for 
simplicity. For each Beltrami
differential $\mu \in B_1 (G_\varphi)$ for 
$G_\varphi$, we pull back
$\mu$ via the developing map $W_\varphi$ to a Beltrami 
differential $\hat
\mu$, namely, we set

$$ \hat \mu = \mu \circ  W_\varphi \times 
 \overline {W_\varphi'}/ W_\varphi'.
$$

The differential $\hat \mu$ is a $(-1, 1)$-form of 
the local
homeomorphism
$f^\mu \circ W_\varphi$, where $f^\mu$ is a 
quasiconformal mapping with
Beltrami differential $\mu$. It is easy to see that 
$\hat \mu$ is an
element of $B_1 (\Gamma, L)$, and the assignment $\hat 
\Phi :   \mu
\to  \hat \mu $ is a holomorphic mapping from 
$B_1(G_\varphi)$ to
$B_1 (\Gamma, L)$.

\proclaim {Proposition 3.2} In the above situation, 
the mapping
$\hat \Phi : \mu \to \hat \mu $ descends to a 
holomorphic mapping
$\Phi : QC(G_\varphi) \to T(\Gamma)$. (Here, we regard 
$T(\Gamma)$ as a
complex manifold.)
\endproclaim 
\demo
{ Proof} We have to show that two 
equivalent elements
$\mu_1$ and
$\mu_2$ in $B_1(G_\varphi)$ are mapped by $\hat \Phi$ 
to equivalent
elements $\hat \mu_1$ and $\hat \mu_2$ in $B_1(\Gamma, 
L)$. \par We first
show that every trivial differential $\tau$ in
$B_1(G_\varphi)$ is pulled back to a trivial 
differential $\hat
\tau$ in $B_1(\Gamma, L)$. Choose three points in
$\Lambda(G_\varphi)$. For each $\mu \in  
B_1(G_\varphi)$, let
$f^\mu$ denote the quasiconformal mapping with 
Beltrami differential
$\mu$ which fixes the three points. If $\tau$ is 
trivial, then
$f^\tau$ fixes every point on $\Lambda(G_\varphi)$ and 
maps each
component of
$\Omega(G_\varphi)$ onto itself. 
 Here we note how the locally univalent mapping 
$W_\varphi$ for
$\varphi \in \operatorname {int}_{Q(\Gamma)} K(\Gamma)$ behaves in 
$L$. The preimages
$(W_\varphi)^{-1}(\Omega (G_\varphi))$ and
$(W_\varphi)^{-1}(\Lambda (G_\varphi))$ are 
$\Gamma$ invariant
and mutually disjoint.
 By Lemma 3.1, the holonomy representation 
$G_\varphi$
is either a quasifuchsian group or a totally 
degenerate group. In either
case, each component of $\Omega 
(G_\varphi)$ is simply
connected. Hence  each germ of the local inverse of 
$W_\varphi$ has an
analytic continuation in each component of $\Omega
(G_\varphi)$. Therefore each component of 
the preimage 
$(W_\varphi)^{-1}(\Omega (G_\varphi))$ is a 
simply connected
domain which is mapped onto a component of
$\Omega (G_\varphi)$ by $W_\varphi$ 
injectively. 
Decompose $W_\varphi ^{-1}(\Omega(G_\varphi))$ to
  the disjoint union of  connected components and denote it by $\cup_{n = 1}
^{\infty} D_{\varphi,n}$. As we have seen above for each $n$ there is a
branch
$W_{\varphi,n} ^{-1}$ of $W_\varphi ^{-1}$ on $W_\varphi (D_{\varphi,n})$
such that $W_{\varphi,n} ^{-1} \circ W_\varphi = id.$ on $D_{\varphi, n}$.

We define a 
mapping $f: L\to L$ by
$$   f(z)  = \cases z &\quad \text{on $
                  W_\varphi ^{-1} (\Lambda 
(G_\varphi))$} \\  W_{\varphi,n} ^{-1} \circ f^\tau \circ W_\varphi (z)
&\quad \text{on $D_{\varphi, n} \quad (n = 1,2,...)$}. \endcases
$$

Now we claim
that $f$ is a quasiconformal homeomorphism of $L$. 
(Intuitively, $f$ is
nothing but
$f^\tau$ with the local charts given by $W_\varphi$ 
rather than the
standard coordinate of the Riemann sphere.) We have to 
care about the
behavior of $f$ near $W_\varphi ^{-1} (\Lambda 
(G_\varphi))$. \par
Fix a
point $z_0 \in W_\varphi ^{-1} (\Lambda (G_\varphi))$ arbitrarily
and take a small disc
$\delta$ around $z_0$. The image $W_\varphi (\delta)$ 
is a neighborhood
of $W_\varphi (z_0) \in  \Lambda (G_\varphi)$.  Take a 
disc
$\delta'\ni z_0$ small enough so that 
$f^\tau (W_\varphi (\delta')) \subset W_\varphi 
(\delta)$. Recall that 
$f^\tau | W_\varphi (\delta) \cap \Lambda (G_\varphi) 
=  id.$ We can take
the branch $W_{\varphi,\delta'}^{-1}$ of $W_\varphi^{-1}$ on $W_\varphi
(\delta')$ such that
$W_{\varphi,\delta'} ^{-1} \circ
W_\varphi $ is the identity on $\delta'$. With this 
branch of $W_\varphi
^{-1}$, we can define a quasiconformal mapping by 
$W_{\varphi,\delta'} ^{-1}\circ f^\tau \circ W_\varphi$ on 
$\delta'$. \par
Obviously, $W_{\varphi,\delta'} ^{-1}\circ f^\tau \circ W_\varphi$ 
coincides with $f$ on $W_\varphi ^{-1}(\Lambda (G_\varphi))$. We claim
that this mapping coincides with the mapping
$W_{\varphi,n} ^{-1} \circ f^\tau \circ W_\varphi$ also on each component
of 
$D_{\varphi,n} \cap  \delta'$. Then it will follow that $f 
=  W_{\varphi,\delta'} ^{-1}\circ f^\tau \circ W_\varphi$  is
quasiconformal near $z_0$.  (Note that this claim is not  quite obvious,
because 
$W_{\varphi}(\delta') \cap \Omega (G_\varphi)$ may have infinitely many 
components and $f^\tau$ might exchange these components.)
 In order to show the above claim, it suffices to show 
that
$f^\tau$ maps each component of
$W_\varphi (\delta')
\cap
\Omega (G_\varphi)$ into the component of $W_\varphi 
(\delta) \cap \Omega
(G_\varphi)$ containing it. To verify this, we take an 
isotopy $\{f_t\}_{
0\le t
\le 1 }$ between $f^\tau$ and the identity map via 
quasiconformal
mappings fixing every point of
$\Lambda(G_\varphi)$, here, the existence of such an 
isotopy is due to 
Earle-McMullen [EM]. Replacing $\delta'$ into a 
smaller disc, if
necessary, we may assume that
$f_t(W_\varphi(\delta')) \subset  W_\varphi(\delta)$ 
for all $t \in
[0,1]$. Then for each component $u$ of 
$W_\varphi(\delta') \cap \Omega
(G_\varphi)$ and $t$, 
$f_t(u)$ is contained in a 
component of  $W_\varphi(\delta ) 
\cap \Omega
(G_\varphi)$ and the component
 depends continuously 
on $t$. Hence
 it is independent of $t$. \par Thus we have seen that 
$f$ is a
quasiconformal homeomorphism of $L$ with Beltrami 
differential $\hat
\tau$. Since $f$ is the identity map on 
$W_\varphi ^{-1} (\Lambda (G_\varphi))$ and since this 
set is
$\Gamma$-invariant, the extension of $f$ to the real 
axis is the
identity. Hence $\hat \tau$ is a trivial differential. 
\par Now we
proceed to the general case. Assume that $\mu_1$ and 
$\mu_2$ are
equivalent in $B_1 (G_\varphi)$. Then there is a 
trivial differential
$\tau$ in $B_1 (G_\varphi)$ such that 
$f^{\mu_2}=f^{\mu_1}\circ f^{\tau}$. Now, as we have 
shown above,
$\hat \tau $ is trivial in $B_1 (\Gamma, L)$ and that
$W_\varphi \circ f_{\hat \tau } = f^\tau \circ 
W_\varphi$,  where
$f_{\hat \tau }$ is the quasiconformal homeomorphism 
of $L$ with Beltrami
differential $\hat \tau$ fixing $0,1$ and $\infty$. 
Hence we have
$$ \align f^{\mu_2}\circ W_\varphi & = f^{\mu_1} \circ 
f^\tau \circ
W_\varphi\\ & = f^{\mu_1}\circ W_\varphi \circ f_{\hat 
\tau}.
\endalign
$$ Let $f_{\hat \mu _i}$, $i = 1,2$, denote the 
quasiconformal
homeomorphism of $L$ with Beltrami differential $\hat 
\mu _i$ fixing
$0,1$ and $\infty$. Since the Beltrami differential 
$\hat \mu_i$ is the
Beltrami differential of $f^{\mu_i}\circ W_\varphi$  
by definition, the
above equality implies
$$f_{\hat \mu _2} = f_{\hat \mu _1}\circ f_{\hat 
\tau}.$$ Since $\hat
\tau$ is trivial, ${\hat \mu _1}$ and ${\hat \mu _2}$ 
are equivalent.
\qed

\enddemo

 We define another mapping which is essentially 
obtained in the above
argument. 
\par 
Fix a point $\varphi \in 
\operatorname{int}_{Q(\Gamma)} K(\Gamma)$.
For each $\mu \in B_1(G_\varphi)$, let $f_{\hat \mu}$ 
denote the
quasiconformal homeomorphism of $L$ with Beltrami 
differential fixing
$0,1$ and $\infty$.  Then by the definition of  ${\hat 
\mu}$, the mapping
$W^\mu$ defined by $W^\mu = f^\mu \circ W_\varphi \circ 
f_{\hat \mu}^{-1}$ is a
locally univalent holomorphic mapping on $L$ to the 
Riemann sphere.
Moreover,
$$ W ^\mu \circ (f_{\hat\mu} \circ \gamma \circ f_{\hat 
\mu}^{-1})
 = (f^\mu\circ \chi_\varphi(\gamma) \circ (f^\mu)^{-
1})\circ W^\mu,\tag 3.1$$ 
for all $\gamma \in \Gamma$. Thus $W^\mu$ is a developing map of a
projective structure on $L/\Gamma_{\hat \mu}$ where $\Gamma _{\hat \mu} = 
f_{\hat  \mu}\Gamma f_{\hat \mu}^{-1}$.  

Taking the Schwarzian derivative of $W^\mu$,
we have a holomorphic quadratic differential 
$\tilde\phi (\mu) \in Q(\Gamma _\mu)$. 
Thus we have a mapping $\tilde\phi :B_1(G_\varphi) \to Q$. 
 The
mapping
$\tilde\phi: B_1 (G_\varphi) \to Q$  
descends to
a mapping from
$QC(G_\varphi)$ into $Q$. In fact, if $\mu_1$ and $\mu_2$ are
equivalent, $\hat \mu_1$ and $\hat \mu_2$ are equivalent by Proposition
3.2. Therefore we have $\Gamma_{\hat \mu_1} = \Gamma_{\hat \mu_2}$ and
$\tilde \phi (\hat \mu _1)$ and  $\tilde \phi (\hat \mu _2)$ are in the
same fiber $Q(\Gamma_{\hat \mu_1})$. Furthermore, these two quadratic
differentials determine projective structures with the same holonomy
representation. Hence
by the  following Theorem
by Poincar\'e (cf. Kra [Kr1]), these holomorphic 
differential are the same.
\proclaim {Theorem (Poincar\'e)}  Two different 
projective structures on
the same complex structure of a compact surface have 
different holonomy
representations.
\endproclaim 
We shall denote the mapping $QC(G_{\varphi}) \to 
Q$ by 
$\phi$. The mapping
$\phi :QC(G_\varphi) \to Q$ is a continuous mapping which satisfies
$\phi \circ h = id.$, where $h$ is the holonomy map. 
By Theorem 2.5, $\phi$ is a $C^1$ map. 
 Also,
it satisfies 
$\Phi =\pi\circ\phi$, where 
$\pi : Q \to T(\Gamma)$ denotes the projection and $ 
\Phi : QC(G_\varphi) \to T(\Gamma)$ is the holomorphic mapping  defined in 
Proposition 3.2. \par

\proclaim {Proposition 3.3} Let $\varphi$ be an element
of
$\operatorname{int}_{Q(\Gamma)} K(\Gamma)$ and  
let $\kappa_\varphi (\Gamma)$ denote the component of
$\operatorname{int}_{Q(\Gamma)} K(\Gamma)$ containing $\varphi$. Then
$\kappa_\varphi (\Gamma)$ is a
subset of a component of 
$\phi(QC(G_\varphi)) \cap Q(\Gamma)$.
 
\endproclaim 
\demo {Proof} 

We shall show every element $\psi \in  \kappa_\varphi (\Gamma)$ is
contained in $\phi(QC(G_\varphi)) \cap Q(\Gamma)$.   
First, we consider
the case  that
$\psi$ is close to $\varphi$ so that  there is  a 
complex analytic disc
connecting them, namely, there is a family 
$\{\varphi_z\}_{z \in \Delta}
\subset
\kappa_\varphi(\Gamma)$ depending holomorphically on 
$z$ such that
$\varphi = \varphi_0$ and $\varphi_\zeta = \psi$ for 
some $\zeta \in
\Delta$ (Here, $\Delta$ denotes the unit disc of the 
complex plane). Then
the set of holonomy representations
$\chi_{\varphi_z}$ depends holomorphically on $z$. Note that
$\chi_{\varphi_z}$ is isomorphic by Lemma 3.1. 
Applying Theorem 2.9  to the family $\chi_{\varphi_z} \circ
\chi_{\varphi_0}^{-1}$, we have a family of Beltrami differentials 
$\{\mu_z\}$  depending
holomorphically on $z$ such that
$$\chi_{\varphi_z} (\gamma) = f^{\mu_z} \circ 
\chi_\varphi (\gamma)
\circ (f^{\mu_z})^{-1},
$$ for all $\gamma \in \Gamma$ and $z\in \Delta$. For 
the Beltrami 
differential $\hat \mu_z = \hat \Phi (\mu_z)$ and the 
quadratic
differential $\phi (\mu_z) \in Q(\Gamma_{\hat \mu_z})$ we 
shall show that
$\phi (\mu_z) = \varphi_z$ for all $z\in \Delta$. The 
set
$\{z \in \Delta; \phi (\mu_z) = \varphi_z\}$ is 
clearly a closed set  and
non-empty. (It contains $0$.) On the other hand, since 
$\phi (\mu_z)$ and $\varphi_z$ have the same holonomy representation and
the representation
map $h : Q \to Rep$ is a local homeomorphism (Theorem 2.5), the 
set $\{z \in \Delta;
\phi (\mu_z) = \varphi_z\}$ is also an open set. Hence 
$\phi (\mu_z) =
\varphi_z$ for all $z \in
\Delta$. \par
Thus we have shown that  all $\psi \in
\kappa_\varphi(\Gamma)$  sufficiently close
to $\varphi$ is contained  in $Q(\Gamma) \cap \phi(QC(G_\varphi))$. For
general
$\psi\in
\kappa_\varphi(\Gamma)$, we take a chain of 
holomorphic discs connecting
$\psi$ and $\varphi$ and apply the above 
argument to each disc. 
\qed
\enddemo

 Now we can 
prove the main
theorem.
\proclaim {Theorem 3.4} Every component of 
$\operatorname{int}_{Q(\Gamma)} K(\Gamma)$
consists of quasifuchsian groups. Every such component 
is holomorphically
equivalent to a complex analytic submanifold of $T_g \times T_g$.
\endproclaim 
\demo{Proof } To prove the first claim, it is 
sufficient to show that
there is no component of $\operatorname{int}_{Q(\Gamma)} K(\Gamma)$
consisting of totally degenerate 
groups, by Lemma 3.1.
 \par
 Assume that there is a component $\kappa(\Gamma)$ in
$\operatorname{int}_{Q(\Gamma)} K (\Gamma)$ which consists of 
totally degenerate
groups. We will draw a contradiction.  Take an
element 
$\varphi \in \kappa(\Gamma)$ and denote the holonomy group by 
$G_\varphi$ for short. We define a
mapping $\Phi : QC(G_\varphi) \to T(\Gamma)$ as in 
Proposition 3.2.  Then
since
$\kappa(\Gamma)$ and
$QC(G_\varphi)$ are of the same dimension, 
$h(\kappa(\Gamma))$ is  an open set of
$ QC(G_\varphi) $.  
It follows that the holomorphic
mapping $\Phi : QC(G_\varphi) \to T(\Gamma)$ takes 
the value
$\Gamma$ on the open set  
$h (\kappa(\Gamma))$ of $QC(G_\varphi)$.
Therefore it is a constant mapping:
$\Phi (QC(G_\varphi) ) \equiv \Gamma$. It follows that  
$\phi (QC(G_\varphi) ) \subset Q(\Gamma)$.  As $QC(G_\varphi)$ and
$Q(\Gamma)$ are of the same dimension,  $\phi (QC(G_\varphi))$ is an open
subset of $Q(\Gamma)$.  From this fact and Proposition 3.3, we have $\phi
(QC(G_\varphi)) = \kappa (\Gamma)$.
 \par
 It follows that there is a point $\psi \in \kappa(\Gamma)$ 
such that the
holonomy representation $\chi_\psi (\Gamma) \in 
QC(G_\varphi)$ has 
Fuchsian equivalent $\Gamma$. Namely, we can choose a 
point 
$f^\mu  G_\varphi (f^\mu)^{-1} \in QC(G_\varphi)$ such 
that there is a
univalent holomorphic mapping $W : L \to 
\Omega (f^\mu  G_\varphi (f^\mu)^{-1})$ such that $W\circ \gamma = f^\mu  
\chi_\varphi(\gamma)
(f^\mu)^{-1} \circ W$ for all $\gamma \in \Gamma$  and
  $\psi = \phi (\mu)$. Now, the projective structure 
defined by the
univalent map $W$ and the projective structure defined 
by $W_\psi$ have the
same holonomy representation and have the same 
underlying complex
structure. By Poincar\'e's theorem, $W$ and $W_\psi$ coincides up to
M\"obius transformation. However, this is impossible because
$W(L) = \Omega (f^\mu  G_\varphi
(f^\mu)^{-1}) \subsetneqq \hat{\operatorname {\bold C}}$, while
$W_\psi (L) = \hat{\operatorname {\bold C}}$ by Proposition 2.7.
(Note that $\kappa (\Gamma)$ is not the Bers
slice.)  
\par

We have seen each component of $\operatorname{int}_{Q(\Gamma)} K(\Gamma)$
consists of projective structures whose holonomy representations are
quasifuchsian.  From this fact, we shall show that $\kappa _\varphi
(\Gamma)$ coincides with the component of
$\phi(QC(G_\varphi)) \cap Q(\Gamma)$ containing $\varphi$.
Recall the proof of Proposition 3.3. For each
point $\varphi \in \operatorname{int}_{Q(\Gamma)} K(\Gamma)$, the
component
$\kappa_\varphi (\Gamma)$ of $\operatorname{int}_{Q(\Gamma)} K(\Gamma)$
containing $\varphi$ is an
open subset of a component of
$\phi(QC(G_\varphi))  \cap Q(\Gamma)$. On the other hand, as $G_\varphi$
is a quasifuchsian group, $Q$ and $QC(G_\varphi)$ is the same dimension.
Therefore, $h (QC(G_\varphi))$ is an open subset of $Q$, hence
each component of $\phi(QC(G_\phi)) \cap Q(\Gamma)$
is an open subset of $\operatorname{int}_{Q(\Gamma)} K(\Gamma)$.
Therefore, $\kappa _\varphi (\Gamma)$ coincides with the component of
$\phi(QC(G_\varphi)) \cap Q(\Gamma)$ containing $\varphi$. \par

It follows that $\kappa_\varphi
(\Gamma)$ is biholomorphically equivalent to the component of
$\Phi ^{-1} (\Gamma)$ containing $G_\varphi$, which is a submanifold of
$QC(G_\varphi)$ by Lemma 2.6.   Now the second claim of the theorem
follows from the fact that
$QC(G_\varphi)$ is holomorphically equivalent to $T_g \times T_g$.
\qed
\enddemo

 What we have done above is to embed 
$QC(G_\varphi)$ into
$Q$ with the mapping $\phi$ (in fact, we have seen 
that 
$\phi$ is a local inverse to the holonomy map $h : Q 
\to Rep$.) and
consider the slice along  a fiber  of 
$Q \to T_g$. It is easy to see that for every $\Gamma'
\in T(\Gamma)$ each component of the slice over 
$\Gamma'$ is  a component
of $\operatorname{int}_{Q(\Gamma)} K(\Gamma')$ \par

\proclaim {Corollary 3.5} Each component of 
$\operatorname{int}_{Q(\Gamma)} K(\Gamma)$
is complete with respect to the Kobayashi
\hyphenation{Ko-ba-ya-shi} hyperbolic metric and Carath\'eodory metric.
Hence each such component  is a domain of
holomorphy.
\endproclaim 
\demo { Proof} We showed in the above 
theorem that each
component
$\kappa(\Gamma)$ of $\operatorname{int}_{Q(\Gamma)} K(\Gamma)$ is 
holomorphically
equivalent to a  submanifold   $T_g \times T_g$.  It
is well-known that $T_g$  is complete with
respect to the Carath\'eodory metric and Kobayashi 
metric. By a standard
argument  (see Kobayashi [K]),
$T_g \times T_g$ and its submanifolds are complete 
with respect to these
metrics. It is also known that in general a domain in 
$\text{\bf C}^N$ complete with respect to 
Carath\'eodory distance is a
domain of holomorphy.
\qed
\enddemo

\bigskip
Here we give a brief remark on the existence of components of
$\operatorname {int}_{Q(X)} K(X)$ other than the Bers slice, which is
discussed in \cite {T}. 
In the above discussion, we have shown that any component of
$\operatorname{int} _{Q(\Gamma)} K(\Gamma)$ is a component of
$\operatorname{int}  K \cap Q(\Gamma)$. It is shown  in \cite {T} 
that 
$\operatorname{int} K$ has infinitely many component and each component
has non-empty intersection with $Q(\Gamma)$. Therefore we have:
\proclaim {Theorem 3.6} On any complex structure $X 
\in T_g$ there are
infinitely many components of $\operatorname {int}_{Q(X)} 
K(X)$.
\endproclaim

\enddocument